\input amstex
\documentstyle{amsppt}
\magnification=\magstep1 \hsize=5.2in \vsize=6.7in

\catcode`\@=11 \loadmathfont{rsfs}
\def\mycal{\mathfont@\rsfs}
\csname rsfs \endcsname

\topmatter
\title  ON OZAWA'S PROPERTY \\
FOR FREE GROUP FACTORS\endtitle
\author SORIN POPA \endauthor

\rightheadtext{Free group factors}

\affil University of California, Los Angeles\endaffil

\address Math.Dept., UCLA, LA, CA 90095-155505\endaddress
\email popa@math.ucla.edu\endemail

\thanks Supported in part by NSF Grant 0601082.\endthanks

\abstract We give a new proof of a result of Ozawa showing that if a
von Neumann subalgebra $Q$ of a free group factor $L\Bbb F_n, 2\leq
n\leq \infty$ has relative commutant diffuse (i.e. without atoms),
then $Q$ is amenable.
\endabstract

\endtopmatter

\document

\heading \S 1. Introduction \endheading

The purpose of this note is to present a rather short and elementary
proof of the following result of Ozawa:

\proclaim{Theorem 1.1 ([O])} Let $Q\subset N=L\Bbb F_n$ be a von
Neumann subalgebra with $Q'\cap N$ diffuse. Then $Q$ is amenable.
\endproclaim

In fact, Ozawa proved the above property for all group factors
$N=L\Gamma$ with $\Gamma$ word-hyperbolic group. The property
implies in particular that all such factors $N$ are {\it prime},
i.e. they cannot be decomposed as a tensor product of II$_1$
factors. Ge had already proved the primeness of $N=L\Bbb F_n$ in
([G]), by using Voiculescu's free entropy techniques ([V]). A
combination of results in ([GSh]) and ([J1]) showed that in fact all
free products of diffuse von Neumann algebras embeddable into
$R^\omega$ are prime (see also [J2]). Recently, Peterson used a
completely new approach to prove the primeness of group factors
$N=L\Gamma$ with $\Gamma$ satisfying $\beta^{(2)}_1(\Gamma)\neq 0$
([Pe]). His results cover also the primeness of all free products of
diffuse finite von Neumann algebras and of all non-amenable
subfactors of $L\Bbb F_n$.

Both Ozawa's insight, the use of free entropy yechniques in ([G],
[GSh], [J1,2]) and Peterson's new approach are conceptually and
technically involved, and aim at showing the property for as many
factors as possible.

In turn, our proof here targets the free group factor case only. The
argument is based on the following two properties of $N=L\Bbb F_n$:

\vskip .05in (1) A ``malleability deformation'' which allows moving
$N=N*\Bbb C$ along a path $\alpha$ of automorphisms inside $N*N$,
from position $N=N*\Bbb C$ to $\Bbb C * N$. It is important for us
here that the path $\alpha$ be ``group-like'' and to have a certain
``symmetry''. \vskip .05in (2) The fact that any $Q\subset N$ which
has no amenable summand has ``spectral gap'' with respect to the
orthogonal of $N$ in $N*N$. In other words, there exists a
``critical'' finite set of unitaries $F\subset Q$ such that if $x\in
N*N$ almost commutes with all $u\in F$ then $x$ is almost contained
in $N$. \vskip .05in Property (1) was already considered in ([P1])
and, exactly in the form we need here, on (page 322 of [P2]). A
related deformation for free products of arbitrary factors has been
used extensively in ([IPeP]). Property (2) is quite evident by one
of Connes' criteria for amenability, since $L^2(N*N)\ominus L^2 N$
is a multiple of the $N$ bimodule $L^2 N \overline{\otimes} L^2 N$,
of Hilbert-Shmidt operators on $L^2 N$. It may in fact have been
noticed before, in some equivalent form. We will nevertheless prove
it, for completeness, and in fact also include the proof of (1) from
([P2]).

The proof of Theorem 1.1 then goes as follows: By contradiction we
may assume $Q$ has no amenable summand, so by (2) a critical set
$F\subset Q$ can be chosen. Thus, any element in $N*N$ that
approximately commutes with $F$ in $N*N$ must be almost contained in
$N=N*\Bbb C$. Due to (pointwise) continuity of $\alpha_t$, if $t>0$
is small then $\alpha_t(u)$ is close to $u$, $\forall u\in F$, so
elements that commute with $\alpha_t(F)$ must almost commute with
$F$, thus have to be almost in $N$. This shows that the unit balls
of $\alpha_t(Q'\cap N)$ and $Q'\cap N$ are uniformly close one to
another. Known results from ([P1,3] or [PSiSm]) then imply that
$Q'\cap N$ and $\alpha_t(Q'\cap N)$ can be ``intertwined'' with a
non-zero partial isometry. The group-like structure of $\alpha$ and
its symmetry make possible to ``patch'' the intertwiners between
position $t$ to position $t+\Delta t$, for ``incremental'' $\Delta
t$, into a non-zero intertwiner between a diffuse von Neumann
subalgebra of $N * \Bbb C$ and its image under $\alpha_1$ in $\Bbb C
* N$. This in turn contradicts (Corollary 4.3 in [P4]).

Note that this line of proof is identical to the approach in (4.3.2
of [P1]) to Connes-Jones theorem on non-embeddability of property
(T) II$_1$ factors $Q_0$ into free group factors $N$. However, that
argument showed only that there exist no embeddings $Q_0\subset N$
with $Q_0'\cap N=\Bbb C$, and could not settle the case $Q_0'\cap N$
diffuse, because of a poorer technique for handling the ``patching''
of the incremental intertwiners in this case: if $Q_0'\cap N$ is
diffuse and $\alpha_t(Q_0'\cap N)$ is uniformly close to $\alpha_{t
+ \Delta t}(Q_0'\cap N)$ then the ensuing intertwiner $v_t$ between
these two algebras is a non-zero partial isometry (not necessarily a
unitary) and multiplying the $v_t$'s for $t=k 2^{-n}$, $\Delta
t=2^{-n}$, may end up giving the 0 intertwiner between $(Q_0'\cap
N)*\Bbb C$ and $\alpha_1((Q_0'\cap N)*1)=\Bbb C * (Q_0'\cap N)$. But
this issue was resolved in ([P2,3]) by the discovery of ``symmetric
paths'', which allows patching intertwiners in a way that makes each
incremental intertwiner have range matching the domain of the next
incremental intertwiner.

Thus, based on deformation/intertwining techniques in ([P2,3]), the
only conceptual novelty in this paper is the idea of using
``spectral gap rigidity'' of embeddings $Q\subset N\subset M=N*N$,
for non-amenable subalgebras $Q$, which allows to argue that any
deformation by automorphisms of $M$ is uniform on $Q'\cap N$. This
idea was first used in ([P5]), but for deformations of the identity
in McDuff factors $M=Q\overline{\otimes} R$ with $Q$ non-$(\Gamma)$,
to prove that such  factors have ``unique McDuff'' decomposition, up
to unitary conjugacy. The proof of that result, which inspired the
approach here, is presented at the end of the paper.

\heading \S 2. Two lemmas \endheading

Before giving the rigorous details of the argument outlined above,
let us first prove properties (1) and (2). (N.B. We will not need
the last part of Lemma 1 below.)

\proclaim{Lemma 2.1 ([P2])} The free group factors $N=L\Bbb F_n,
2\leq n \leq \infty$ have the following property: There exist a
continuous action $\alpha : \Bbb R \rightarrow {\text{\rm
Aut}}(N*N)$ and a period 2 automorphism $\beta \in {\text{\rm
Aut}}(N*N)$ such that

$(1.1)$. $\alpha_1 (N*\Bbb C) = \Bbb C * N$.

$(1.2)$. $N*\Bbb C \subset (N*N)^\beta$.

$(1.3)$. $\beta \alpha_t \beta = \alpha_{-t}, \forall t$.

Moreover, $\alpha, \beta$ can be chosen so that to commute with all
automorphisms of $N*N$ implemented by permuting the same way,
simultaneously, the left and right generators of $\Bbb F_n * \Bbb
F_n$.
\endproclaim
\noindent {\it Proof} To prove the statement it is clearly
sufficient to construct the continuous action $\alpha$ of $\Bbb R$
on $M=N*N$ with the period 2-automorphism $\beta$ of $M$ such that
$\beta \alpha_t \beta = =\alpha_{-t}, \forall t\in \Bbb R$ and such
that there exist isomorphic copies $N_1, N_2$ of $N$ in $M$ which
are mutually free, generate $M$ and satisfy $\alpha_1(N_1)=N_2$,
$N_1 \subset M^\beta$.

To this end, let $a_1, a_2,... $ be the generators of $\Bbb F_n$
viewed as unitary elements in $N*\Bbb C$ and $b_1, b_2,... $ the
same generators but viewed  as unitary elements in $\Bbb C * N$.

Let $h_k \in \Bbb C*N$ be self-adjoint elements with spectrum in
$[0,2]$ such that $b_k ={\text{\rm exp}}(\pi i h_k), \forall k$. We
then put $\alpha_t (a_k) = {\text{\rm exp}}(\pi i t h_k) a_k$ and
$\alpha_t (b_k) = b_k, k=1,2,..., t\in \Bbb R$. It is trivial to see
that ${\text{\rm exp}}(\pi i t h_k) a_k$ and $b_k$ are mutually free
and generate the same von Neumann algebra as $a_k, b_k$. Thus,
$\alpha_t$ defines an automorphism of $N*N$, $\forall t\in \Bbb R$.
Note also that $\alpha_1(a_k)=b_ka_k=b'_k$ is free with respect to
$a_k$ and they jointly generate the same algebra as $a_k, b_k$ do.
Thus, if we let $N_1=N*\Bbb C$ and $N_2=\alpha_1(N*\Bbb C)$ then
$N_1, N_2$ are free and generate $M$. Moreover, by the definition we
clearly have $\alpha_t \alpha_s = \alpha_{t+s}, \forall t,s \in \Bbb
R$, showing that $\alpha$ is a continuous action.

Define now $\beta (a_k) = a_k$ and $\beta (b_k)= b_k^*, \forall
k=1,2,...$. This clearly defines a period 2 automorphism of $M=N*N$
satisfying $N* \Bbb C \subset (N*N)^\beta=M^\beta$. Moreover
$$
\beta (\alpha_t (\beta(a_k))) = \beta (\alpha_t(a_k)) =
\beta({\text{\rm exp}}(\pi i t h_k) a_k)
$$
$$
= {\text{\rm exp}}(\pi i t h_k)^* a_k ={\text{\rm exp}}(-\pi i t
h_k) a_k = \alpha_{-t} (a_k).
$$
Similarly, we get
$$
\beta (\alpha_t (\beta(b_k))) = b_k = \alpha_{-t}(b_k),
$$
showing that all conditions are satisfied.

\hfill $\square$

From the next lemma we only need $\Rightarrow$ in the case
$N=P=L\Bbb F_n$.

\proclaim{Lemma 2.2} Let $Q$ be a separable von Neumann subalgebra
of a finite von Neumann algebra $(N,\tau)$ and $(P,\tau)$ a von
Neumann algebra $\neq \Bbb C$. Then $Q$ has no amenable direct
summand if and only if, when viewing $Q=Q*\Bbb C$, $N=N*\Bbb C$ as
subalgebras of $M=N*P$, we have $Q' \cap M^\omega\subset N^\omega$.
\endproclaim
\noindent {\it Proof}. If $Q$ has amenable direct summand $Qp$, for
some non-zero $p\in \Cal P(Q)$, then by Connes theorem $Qp$ is
approximately finite dimensional and the statement follows from the
case when $Qp$ is finite dimensional, when we trivially have
$(Qp)'\cap pM^\omega p \nsubseteq N^\omega$.

If $Q$ doesn't have amenable direct summand but  $Q'\cap M^\omega
\nsubseteq N^\omega$, then let $x=(x_n)_n\in Q'\cap M^\omega$, $x
\not\in N^\omega$. By substracting $E_{N^\omega}(x)=(E_N(x_n))_n$
(which still commutes with $Q$), we may assume $x \neq 0$ but
$E_N(x_n)=0, \forall n$, thus $x_n\in L^2 M \ominus L^2 N$. By
replacing $(x_n)_n$ with a subsequence and using the separability of
$Q$, we may also assume $[x_n, y]\rightarrow 0$, $\forall y\in Q$.
Note that if $b$ is the conditional expectation of $xx^*$ onto $Q$,
then $b\in Q'\cap Q=\Cal Z(Q)$ and $\lim_\omega
\tau(yx_nx_n^*)=\tau(yb), \forall y\in Q$, is a normal trace on $Q$.

By the usual decomposition of $L^2M=L^2(N*P)$ we see that as $N$
bimodule $L^2M \ominus L^2 N$ is identical to the Hilbert space
$\Cal H$ obtained as the direct sum of infinitely many copies of
$L^2 N\overline{\otimes} L^2 N$, i.e. $\Cal H=L^2
N\overline{\otimes} L^2 N \overline{\otimes} \ell^2 \Bbb N$. We
regard each summand $L^2 N \overline{\otimes} L^2N$ as
Hilbert-Schmidt operators on $L^2 N$, and view $\Cal H$ as the space
$HS(L^2 N) \overline{\otimes} \ell^2 \Bbb N$ of ``diagonal''
Hilbert-Schmidt operators in $HS(L^2 N \overline{\otimes} \ell^2
\Bbb N),$ in the obvious way, with $\langle \cdot, \cdot
\rangle_{HS}$ the scalar product on $\Cal H$. If we view the
elements $x_n$ of the sequence $(x_n)_n$ as Hilbert-Schmidt
operators $X_n$ in $\Cal H$, then $\lim_n \|uX_nu^*-X_n\|_{HS}=0,$
$\forall u\in \Cal U(Q)$ and $\lim_\omega \langle y X_n,X_n
\rangle_{HS}=\lim_\omega \tau(yx_nx_n^*)=\tau(yb).$ Thus, if we
define $\varphi (T) \overset \text{\rm def} \to = \lim_\omega
\langle T X_n, X_n \rangle_{HS}$, $\forall T\in \Cal B=\Cal B(L^2
N\overline{\otimes} \ell^2 \Bbb N)$ then $\varphi$ is a positive
functional on $\Cal B$ which has $Q$ in its centralizer and
restricted to $Q$ equals $\tau( \cdot b)$, i.e. a normal trace. Thus
$Q$ has an amenable summand, by [C]. \hfill $\square$

\heading \S 3. Proof of the main result \endheading

\noindent {\it Proof of Theorem} 1.1. Assume $Q$ has a non-amenable
direct summand. Since $N=L\Bbb F_n$ is a factor, by replacing if
necessary $Q$ with $m \times m$ matrices over a corner of it, we may
assume $Q$ does not have an amenable direct summand but still
$Q^0=Q'\cap N$ is diffuse. By Lemma 2.2, $\forall \varepsilon > 0$,
$\exists \delta > 0, F \subset \Cal U(Q)$ finite such that if $x\in
(M)_1$ satisfies $\|ux-xu\|_2 \leq \delta$, $\forall u\in F$, then
$\|x-E_N(x)\|_2 \leq \varepsilon$.

Let $(\alpha,\beta)$ be as in Lemma 2.1. By continuity of $\alpha$
there exists $t=2^{-n}$ such that $\|\alpha_{-t}(u)-u\|_2 \leq
\delta/2, \forall u\in F$. It follows that
$$
\|[\alpha_{-t}(x), u]\|_2 \leq 2 \|u-\alpha_{-t}(u)\|_2\leq \delta,
\forall x\in (Q^0)_1, u\in F
$$
implying that $\|\alpha_{-t}(x)-E_N(\alpha_{-t}(x))\|_2 \leq
\varepsilon$, $\forall x\in (Q^0)_1,$ or equivalently
$$
\|x-E_{\alpha_t(N)}(x)\|_2 \leq \varepsilon, \forall x\in (Q^0)_1.
$$
By ([P1], or 2.1 in [P3], or [PSiSm]; a modification of arguments in
[Ch] will also do), there exist projections $q\in Q^0$, $q' \in
(Q^0)'\cap M$, $p \in \alpha_t(N)$, an isomorphism $\rho:
qQ^0q\rightarrow p\alpha_t(N)p$, a projection $p'\in
\rho(qQ^0qq')'\cap pMp$ and a partial isometry $v\in M$ such that
$v^*v=qq', vv^*=pp'$, $vx=\rho(x)v, \forall x\in qQ^0q$, and $v$ at
distance $f(\varepsilon)$ from $1$, where $\lim_{\varepsilon
\rightarrow 0} f(\varepsilon)=0$.

But by (Remark 2) of 6.3 in [P4]), since $Q$ is diffuse (because
non-amenable) $Q'\cap M \subset N$. Thus $q'\subset N$ so altogether
$Q_0=qQ^0qq'$ lies entirely in  $N$.  Similarly, since the image of
$\rho$ is a diffuse algebra and $\alpha_t(N)$ splits off $M$ in a
free product, by ([P4]) again it follows that $p'\in \alpha_t(N)$.
We have thus just shown that there exists a diffuse subalgebra
$Q_0=qQ^0qq'\subset N$ and a partial isometry $v$ in $M$ such that
$v^*v=1_{Q_0}$, $vQ_0v^*\subset \alpha_t(N)$.

With $n$ the fixed integer with $2^{-n}=t$ as before, we now
construct by induction over $k\geq 0$ some partial isometries $v_k
\in M$ and diffuse weakly closed von Neumann subalgebras $Q_k
\subset N=N* \Bbb C$ such that

$$
\tau(v_kv_k^*) = \tau(vv^*), v_k^*v_k = 1_{Q_k}, v_kQ_kv_k^* \subset
\alpha_{1/2^{n-k}}(N) \tag a
$$

Letting $Q_0=Q_0, v_0=v$, we see that the relation holds true for
$k=0$. Assume we have constructed $v_j, Q_j$ for $j=0, 1, \ldots
,k$. By applying the automorphism $\beta$ of Lemma 2.2 to the
inclusion in $(a)$ and using the properties of $\beta$, it follows
that
$$
\beta(v_k)Q_k\beta(v_k)^* \subset \alpha_{-1/2^{n-k}}(N) \tag b
$$
By further applying $\alpha_{1/2^{n-k}}$ to this latter inclusion it
follows that
$$
Q_{k+1}\overset\text{\rm def} \to = \alpha_{1/2^{n-k}}(\beta(v_k))
\alpha_{1/2^{n-k}}(Q_k) \alpha_{1/2^{n-k}}(\beta(v_k^*)) \subset
N\otimes \Bbb C \tag c
$$
By conjugating $(c)$ with $\alpha_{1/2^{n-k}}(\beta(v_k^*))$ we thus
get:
$$
\alpha_{1/2^{n-k}}(\beta(v_k^*)) Q_{k+1}
\alpha_{1/2^{n-k}}(\beta(v_k)) = \alpha_{1/2^{n-k}}(Q_k) \subset
\alpha_{1/2^{n-k}}(N) \tag d
$$
On the other hand, by applying $\alpha_{1/2^{n-k}}$ to $(a)$ we also
have:
$$
\alpha_{1/2^{n-k}}(v_k) \alpha_{1/2^{n-k}}(Q_k)
\alpha_{1/2^{n-k}}(v_k^*) \subset \alpha_{1/2^{n-k-1}}(N) \tag e
$$
Altogether, it follows that if we let $v_{k+1} =
\alpha_{1/2^{n-k}}(v_k)\alpha_{1/2^{n-k}}(\beta(v_k^*))$ then by
$(d)$ and $(e)$ we get:
$$
v_{k+1}Q_{k+1}v_{k+1} \subset \alpha_{1/2^{n-k-1}}(N)
$$
Moreover, since $\beta(v_k^*v_k)=v_k^*v_k$, we also have
$v_{k+1}v_{k+1}^* = \alpha_{1/2^{n-k}}(v_kv_k^*)$, so that
$\tau(v_{k+1}v_{k+1}^*) = \tau (v_kv_k^*)$. This ends the induction
argument. Taking $k=n$, by $(1.1)$ it follows that $v_nQ_nv_n^*
\subset \alpha_1(N*\Bbb C) = \Bbb C * N$. But by (Corollary 4.3 in
[P4]), this implies that $v_n=0$. Since $\tau(v_nv_n^*) =
\tau(vv^*)$, it follows that $v=0$, a contradiction. \hfill
$\square$

\heading \S 4. Some remarks\endheading

$1^\circ$. In ([O]) Ozawa called {\it solid} a II$_1$ factor $N$
with all subalgebras with diffuse relative commutant being amenable.
As we mentioned in the introduction, this property implies $N$ is
prime. Primeness phenomena for II$_1$ factors first appeared in
([P4]), where it was shown that free group factors with uncountably
many generators have this property. The natural problem of showing
that separable free group factors are prime as well ([P4]),
mentioned also in ([P6]), remained open for some time, until Ge
solved it using Voiculescu's free entropy techniques ([Ge]). Ge then
asked whether there exist diffuse subalgebras of $L\Bbb F_n$ with
non-amenable relative commutant, a problem settled by Ozawa's
result.

$2^\circ$. Despite strong interest and much effort, Voiculescu's
other remarkable application of free entropy to free group factors,
namely the non-existence of Cartan subalgebras in $L\Bbb F_n$, could
not be shown (so far) by the new methods of Ozawa and Peterson ([O],
[Pe]), nor by Connes-Shlyakhtenko's cohomological invariants ([CS]).
Our approach here may open up some new perspectives in this respect.

$3^\circ$. Ozawa's property can be viewed as a ``non-embeddability
into $L\Bbb F_n$'' of algebras of the form
$\tilde{Q}=\overline{\text{\rm sp}} Q Q'$, with $Q, Q'$ commuting,
$Q$ non-amenable and $Q'$ diffuse. Combining this with an argument
of Connes-Jones in ([CJ]) can be used to show that any group
$\Lambda$ which can be decomposed as a product of two commuting
subgroups $\Lambda = HH'$ with $H$ non-amenable and $H'$ infinite,
has a properly outer cocyle action on $L\Bbb F_\infty$ which cannot
be inner-perturbed to an actual action (i.e. the 2-cocycle involved
does not vanish). Indeed, because if it could then by ([CJ]) we
would have $L_\mu(\Lambda) \subset L\Bbb F_2$, for some scalar
$2$-cocycle $\mu: \Lambda \times \Lambda \rightarrow \Bbb T$ of the
group $\Lambda$, where $L_\mu(\Lambda)$ denotes as in ([CJ]) the von
Neumann algebra generated by the unitaries $u_g$ acting on
$\xi_h=(\delta_{h,k})_k \in \ell^2\Lambda$ by $u_g(\xi_h)=\mu(g,h)
\xi_{gh}$, for $g,h\in \Lambda$. If $H'$ contains torsion free
elements, then it is easy to see that this entails $L_\mu (H)'\cap
(L\Bbb F_2)^\omega$ has a diffuse part. But then (Proposition 7 in
[O]) implies $L\Bbb F_2$ has a diffuse abelian subalgebra with
non-amenable relative commutant in $L\Bbb F_2$, contradicting
Theorem 1. Similarly if $H'$ contains elements of arbitrary large
period. We are thus reduced to the case when there exists $n< \infty
$ such that $g^n=e$, $\forall g\in H'$. This implies the 2-cocycle
$\mu$ satisfies $(\mu(g,h)/\mu(h,g))^n=1$ for all $g\in H', h\in H$.
Thus, if we denote by $\{u_g\mid g \in \Lambda\}$ the canonical
unitaries in $L_\mu (\Lambda)$, then for each $g\in H'$ the subgroup
$H_g=\{h\in H \mid u_hu_g=u_gu_h\}$ is normal with index $\leq n$ in
$H$, and in fact $H/H_g \hookrightarrow \Bbb Z/n\Bbb Z.$ Applying
this to all $g\in H'$ we get a group morphism $H \rightarrow (\Bbb
Z/n\Bbb Z)^{H'}$. Since the right hand group is abelian and $H$ is
non-amenable, the kernel $H_0\subset H$ must be non-amenable. Thus
$Q=L_\mu(H_0)$ is non-amenable and commutes with the diffuse algebra
$L_\mu(H')$, contradicting Theorem 1.

\heading \S 5. A related result \endheading

As mentioned earlier, we end the paper with the proof of a result in
([P5]), on the unique tensor product decomposition of McDuff II$_1$
factors, which inspired the above proof of Theorem 1.1.

Thus, let $M$ be a factor of the form $M=P \overline{\otimes} R$
where $R$ is a copy of the hyperfinite II$_1$ factor and $P$ is a
non$(\Gamma)$ II$_1$ factor. Note that if $M_n\subset R$ are finite
dimensional such that $\overline{\cup_n M_n}=R$, then the factors
$P_n=P \otimes M_n$ increase to $M$ and thus the expectations
$E_{P_n}$ converge to $id_M$.

Let $Q\subset M$ be another non$(\Gamma)$ II$_1$ subfactor of $M$
such that $M=Q \vee Q'\cap M$. By ([C]) it follows that $Q$ has
``spectral gap'' with respect to $Q'\cap M$, i.e. $\forall
\varepsilon>0$, $\exists u_1,...,u_n\in \Cal U(Q)$ and $\delta >0$
such that if $x\in (M)_1$ satisfies $\|[u_i,x]\|_2 < \delta, \forall
i$, then $\|x-E_{Q'\cap M}(x)\|_2 \leq \varepsilon$.

Since $E_{P_n} \rightarrow id_M$, there exists $n$ such that
$\|E_{P_n}(u_i)-u_i\|_2 < \delta/2, \forall i$. But then for any
$x'\in (P'\cap M)_1$ we have
$$
\|[x',u_i]\|_2 =\|[x', u_i-E_{P_n}(u_i)]\|_2 < \delta,
$$
implying that $\|E_{Q'\cap M}(x')-x'\|_2 \leq \varepsilon$. This
shows that $P_n'\cap M \subset_\varepsilon Q'\cap M=R$, which
implies
$$
Q=(Q'\cap M)'\cap M \subset_{2 \varepsilon} (P_n'\cap M)'\cap M=P_n.
$$

By ([OP]) this implies that there exists $u\in \Cal U(M)$ such that
$uQu^*=P^t, u(Q'\cap M)u^*=R^{1/t}$, after some ``re-scaling'' with
an appropriate $t>0$ of the decomposition $M=P\overline{\otimes} R$,
as in ([OP]). We have thus proved:

\proclaim{Theorem 5.1 ([P5])} If a $\text{\rm II}_1$ factor $M$ has
two decompositions of the form $M=N_1 \overline{\otimes} R_1 = N_2
\overline{\otimes} R_2$, with $N_1, N_2$ non$(\Gamma)$ $\text{\rm
II}_1$ factors and $R_1, R_2\simeq R$, then there exists a unitary
element $u\in M$ such that $uN_1u^* = N_2^t$, $uR_1u^*=R_2^{1/t}$,
for an appropriate $t>0$.
\endproclaim

If one applies this to free group factors $N_1=L\Bbb F_n$, $2 \leq n
< \infty$, and  $N_2 = L\Bbb F_\infty$, it follows that $L\Bbb
F_n\overline{\otimes} R \simeq L\Bbb F_\infty\overline{\otimes} R$
iff $L\Bbb F_n\overline{\otimes} \Cal B(\ell^2 \Bbb N) \simeq L\Bbb
F_\infty \overline{\otimes} \Cal B(\ell^2 \Bbb N)$. By a result of
Dykema ([Dy]) and Radulescu ([R]), this shows the following:

\proclaim{Corollary 5.2} The free group factors $L(\Bbb F_s), 2 \leq
s \leq \infty$,  are all isomorphic iff they are isomorphic after
``stabilizing'' with $R$, i.e. iff $L\Bbb F_s \overline{\otimes} R
\simeq L\Bbb F_\infty\overline{\otimes} R$, $\forall s\geq 2$, and
iff there exists $s \geq 2$ such that $L\Bbb F_s \overline{\otimes}
R \simeq L\Bbb F_\infty\overline{\otimes} R$.
\endproclaim

\head  References\endhead

\item{[Ch]} E. Christensen: {\it Subalgebras of a finite algebra}, Math.
Ann. {\bf 243} (1979), 17-29.

\item{[C]} A. Connes: {\it Classification of injective factors},
Ann. Math., {\bf104} (1976), 73-115.

\item{[CJ]} A. Connes, V.F.R. Jones: {\it Property T
for von Neumann algebras}, Bull. London Math. Soc. {\bf 17} (1985),
57-62.

\item{[CS]} A. Connes, D. Shlyakhtenko: {\it $L^2$-homology for
von Neumann algebras}, J. reine angew. Math. {\bf 586} (2005),
125-168.

\item{[D]} K. Dykema: {\it
Interpolated free group factors}, Duke Math J. {\bf 69} (1993),
97-119.

\item{[G]} L. Ge: {\it Applications of free entropy to finite
von Neumann algebras}, Ann. Math. {\bf 147} (1998), 143-157.

\item{[GSh]} L. Ge, J. Shen: {\it On free entropy dimension of
finite von Neumann algebras}, GAFA {\bf 12} (2002), 546-566.

\item{[IPeP]} A. Ioana, J. Peterson, S. Popa: {\it Amalgamated free
products of w-rigid factors and calculation of their symmetry
groups}, math.OA/0505589, to appear in Acta Math.

\item{[J1]} K. Jung: {\it A free entropy dimension lemma}, math.OA/0207149.

\item{[J2]} K. Jung: {\it Strongly} $1$-{\it bounded von Neumann algebras},
to appear in GAFA, math.OA/ \newline 0510576.

\item{[O]} N. Ozawa: {\it Solid von Neumann algebras},
Acta Math. {\bf 192} (2004), 111-117.

\item{[OP]} N. Ozawa, S. Popa: {\it Some prime factorization results
for type} II$_1$ {\it factors}, Invent Math. {\bf 156} (2004),
223-234.

\item{[Pe]} J. Peterson, {\it A $1$-cohomology characterization of property}
(T) {\it in von Neumann algebras}, math.OA/0409527.

\item{[P1]} S. Popa: {\it Correspondences},
INCREST preprint 1986, unpublished, on the web at
www.math.ucla.edu/~popa/preprints.html

\item{[P2]} S. Popa: {\it Some rigidity results for
non-commutative Bernoulli shifts}, J. Fnal. Analysis {\bf 230}
(2006), 273-328.

\item{[P3]} S. Popa: {\it Strong Rigidity of} II$_1$ {\it Factors
Arising from Malleable Actions of $w$-Rigid Groups} I, Invent. Math.
{\bf 165} (2006), 369-408 (math.OA/0305306).

\item{[P4]} S. Popa: {\it Orthogonal pairs of *-subalgebras in
finite von Neumann algebras}, J. Operator Theory {\bf 9} (1983),
253-268.

\item{[P5]} S. Popa: {\it Deformation and rigidity in the study of} II$_1$
{\it factors}, Mini-Course at College de France, Nov. 2004.

\item{[P6]} S. Popa: {\it
Free independent sequences in type} $II_1$ {\it factors and related
problems}, Asterisque {\bf 232} (1995), 187-202.

\item{[PSiSm]} S. Popa, A. Sinclair, R. Smith: {\it Perturbations
of subalgebras of type} II$_1$ {\it factors}, J. Fnal. Analysis,
{\bf 213} (2004), 346-379.

\item{[R]} F. Radulescu: {\it Random matrices,
amalgamated free products and subfactors of the von Neumann algebra
of a free group}, Invent. Math. {\bf 115} (1994), 347--389.

\item{[V]} D. Voiculescu: {\it The analogues of entropy and
of Fisher's information theory in free probability} II, Invent.
Math. {\bf 118} (1994), 411-440 .

\enddocument